\documentclass[12pt]{article}

\usepackage{latexsym}
\usepackage{amssymb}

\newtheorem{theorem}{Theorem}

\newtheorem{corollary}[theorem]{Corollary}
\newtheorem{definition}[theorem]{Definition}

\newcommand{\proof}{\noindent {\it {Proof. }}}
\newcommand{\fproof}{\noindent {\it {First proof. }}}
\newcommand{\sproof}{\noindent {\it {Second proof. }}}

\newcommand{\prend}{ $\diamondsuit $\hfill \bigskip}

\newcommand\id{\mathop{\rm id}}

\newcommand\nph{\varphi}
\newcommand\cch{{\mathcal {H}}}

\newcommand\ccb{{\mathcal {B}}}
\newcommand\ccc{{\mathcal {C}}}

\newcommand\cck{{\mathcal {K}}}

\newcommand\ccr{{\mathcal {R}}}

\newcommand\ccl{{\mathcal {L}}}

\begin{document}


\title{A characterization of Morita equivalence pairs}

\author{I. G. Todorov}
\date{}
\maketitle

\begin{abstract}
We characterize the pairs of 
operator spaces which occur as pairs of Morita equivalence bimodules between 
non-selfadjoint operator algebras
in terms of the mutual relation between the spaces.
We obtain a characterization of the operator spaces which are 
completely isometrically isomorphic to imprimitivity bimodules between some strongly 
Morita equivalent (in the sense of Rieffel) C*-algebras.
As corollaries, we give representation results for such operator spaces.
\end{abstract}

\section{Introduction and preliminaries}

The notion of Morita equivalence is fundamental in both Algebra and Analysis. 
This notion entered the Theory of Operator Algebras with Rieffel's paper \cite{rief}. 
In that work, Rieffel
defined Morita equivalence for C*- and W*-algebras and proved various results 
concerning these notions. After the formulation of the theory of (abstract) operator spaces 
(see \cite{er} for a survey and references), 
this notion was generalized \cite{bmp} for non-selfadjoint 
operator algebras and obtained an important place in the Theory. 
In this paper we give a characterization of those pairs $(X,Y)$ of 
operator spaces, which occur as pairs of Morita equivalence bimodules between 
non-selfadjoint (approximately unital) operator algebras, 
in terms of the mutual relation between the spaces $X$ and $Y$.
We specialize our result to obtain a characterization of the operator spaces which are 
completely isometrically isomorphic to imprimitivity bimodules between strongly 
Morita equivalent (in the sense of Rieffel \cite{rief}) C*-algebras.
We call these (abstract) operator spaces {\it ternary operator systems}. 
The Banach space analogue of these spaces, known as {\it ternary C*-rings} have 
been studied by Zettl \cite{z} and others. Our result is an operator space version of 
an implicit characterization of the imprimitivity bimodules between C*-algebras up to a 
Banach space isometry, contained in \cite{z}.
According to our result, additional conditions must be imposed on an operator space 
which is also a ternary C*-ring, in order to be {\it completely isometric} to a Morita
equivalence bimodule between some C*-algebras.
Thus, the similarity beween ternary C*-rings and ternary operator systems is the 
associativity property of the triple product; the difference is that the \lq\lq C*-condition'' 
for a ternary operator system takes into account its operator space structure as well 
(see Definition \ref{d_tro}).

Finally, we point out generalizations of the Representation Theorem of \cite{z} for the case of a
ternary operator system.

In order to state and prove our results, we recall some notions and notation. We have 
as a main reference the monograph \cite{er}. 
An {\it operator space} is a Banach space $X$ endowed with a sequence of norms on
the linear space $M_n(X)$ of $n\times n$ matrices with entries in $X$
for which the (left and right) actions of $M_n(\mathbb{C})$ on $M_n(X)$ are contractive and
the norm of a direct sum of two matrices equals the maximum of their norms. A map 
$\nph : X\longrightarrow Y$ between the operator spaces $X$ and $Y$ 
is said to be {\it completely bounded (contractive, isometric)}, 
if the map $\nph_n : M_n(X)\longrightarrow M_n(Y)$ given by $\nph_n((x_{i j})) = (\nph(x_{i j}))$
(where $(x_{i j})\in M_n(X)$)
is bounded (contractive, isometric) for each $n\in\mathbb{N}$ and 
$\|\nph\|_{cb} = \sup_{n\in\mathbb{N}}\|\nph_n\|$ is finite.
If $x = (x_{i j}), y = (y_{i j})\in M_n(X)$, then 
we let $x\odot y = (\sum_k x_{i k}\otimes y_{k j})_{i, j}\in M_n(X\otimes Y)$. 
If an operator space $X$ is an $A, B$-bimodule, where $A$ and $B$ are some 
algebras, we define 
$a\odot x = (\sum_k a_{i k}\cdot x_{k j})_{i, j}$ and 
$x\odot b = (\sum_k x_{i k}\cdot b_{k j})_{i, j}$, 
where $a = (a_{i j})\in M_n(A)$ and $b = (b_{i j})\in M_n(B)$.
The {\it Haagerup tensor product} 
of the operator spaces $X$ and $Y$ can be described by defining a norm $\|\cdot\|_n$ on 
$M_n(X\otimes Y)$, setting 
$$\|U\|_n = \inf \{\|x\|\|y\| : U = x\odot y, x\in M_{n,p}(X), y\in M_{p,n}(Y), p\in\mathbb{N}\}$$
for each $U\in M_n(X\otimes Y)$ and 
$n\in\mathbb{N}$ and letting $X\otimes_hY$ be the completion of $X\otimes Y$ 
with respct to $\|\cdot\|_1$. The Haagerup tensor product is associative. 
Moreover, if $X_1,X_2,Y_1$ and $Y_2$ are operator spaces and $f_1:X_1\rightarrow Y_1$ and 
$f_2:X_2\longrightarrow Y_2$ are completely bounded maps then there is a unique completely 
bounded map $f_1\otimes f_2 : X_1\otimes_h X_2\longrightarrow Y_1\otimes_h Y_2$ such that 
$f_1\otimes f_2 (x_1\otimes x_2) = f_1(x_1)\otimes f_2(x_2)$, $x_1\in X_1$, $x_2\in X_2$.
We have that $\|f_1\otimes f_2\|_{cb}\leq \|f_1\|_{cb}\|f_2\|_{cb}$.
Each multilinear completely bounded map from the direct product of several operator spaces 
into an operator space can be linearized through the Haagerup tensor product of its domain 
spaces (see \cite{cs} for the theory of multilinear completely bounded maps). 
An approximately unital Banach algebra will be a Banach algebra which is not unital but 
which possesses a contractive approximate identity. We assume that the norm of the unit in a
unital Banach algebra is one.
If an approximately 
unital or unital Banach algebra is an operator space and its multiplication is completely contractive 
(we call such an object an {\it (abstract) operator algebra}),
then it is completely isometrically isomprphic to an algebra of operators on a Hilbert space.
Two approximately unital or unital operator algebras $A$ and $B$ are said to be Morita equivalent if 
there exist operator spaces $X$ and $Y$, such that 1) $X$ is an $A,B$-bimodule, $Y$ is a
$B,A$-bimodule with completely contractive module actions; 2) there are completely bounded
bimodule maps
$(\cdot,\cdot) : X\times Y\longrightarrow A$ and  $[\cdot,\cdot] : Y\times X\longrightarrow B$ 
such that $(x\cdot b,y) = (x, b\cdot y)$, $[y\cdot a,x] = [y, a\cdot x]$ (that is, these maps are 
{\it balanced}), $(x_1,y)\cdot x_2 = x_1\cdot [y,x_2]$,
$[y_1,x]\cdot y_2 = y_1\cdot (x,y_2)$ for each $x,x_1,x_2\in X$, $y,y_1,y_2\in Y$, $a\in A$,
$b\in B$ and 
3) the linearized maps on the Haagerup tensor 
product induced by the pairings are complete quotients (see \cite{bmp} for the exact definition of this
condition). The 6-tuple $(A,B,X,Y,(\cdot,\cdot),[\cdot,\cdot])$ is called a {\it Morita context}.

If $X$ is a Banach space, we denote by $B(X)$ the (Banach) algebra of bounded linear 
operators on $X$.

\section {Morita pairs between non-selfadjoint operator algebras}

In this section we characterize Morita pairs between non-selfadjoint approximately unital 
or unital operator algebras. We give a full proof of the result only in the case both algebras 
are approximately unital. We then comment the necessary changes to be made if one (both)
of the algebras is (are) unital.

\begin{theorem}\label{th_mor_pair}
Let $X$ and $Y$ be operator spaces. Then $(X,Y)$ is a Morita equivalence pair between some 
approximately unital operator algebras 
if and only if there exist completely bounded maps
$p: X\otimes_h Y\otimes_h X\rightarrow X$ and
$q: Y\otimes_h X\otimes_h Y\rightarrow Y$ and nets
$\{e_{\alpha}\}\subset X\otimes Y$ and 
$\{f_{\beta}\}\subset Y\otimes X$ with $\|e_{\alpha}\|_{cb} < 1$ and $\|f_{\beta}\|_{cb} < 1$ 
such that 

1) $p(p(x_1\otimes y_1\otimes x_2)\otimes y_2\otimes x_3)=
p(x_1\otimes q(y_1\otimes x_2\otimes y_2)\otimes x_3)=
p(x_1\otimes y_1\otimes p(x_2\otimes y_2\otimes x_3))$;

2) $q(q(y_1\otimes x_1\otimes y_2)\otimes x_2\otimes y_3)=
q(y_1\otimes p(x_1\otimes y_2\otimes x_2)\otimes y_3)=
q(y_1\otimes x_1\otimes q(y_2\otimes x_2\otimes y_3))$;

3) $\lim_{\alpha}p(e_{\alpha}\otimes x) = x$ and

4) $\lim_{\beta}q(f_{\beta}\otimes y) = y$

\noindent for each $x,x_i\in X$ and $y,y_i\in Y$ ($i=1,2,3$). Moreover, the pairings in the 
respective Morita context are completely contractive if and only if the
mappings $p$ and $q$ with the properties 1)-4) can be chosen to be completely contractive.
\end{theorem}
We will give two proofs of the above theorem. The first one is 
applicable only for completely contractive pairings, while the second one can be applied 
in the general case.

\bigskip

\fproof
Let $X$ and $Y$ be operator spaces, 
$p: X\otimes_h Y\otimes_h X\rightarrow X$ and
$q: Y\otimes_h X\otimes_h Y\rightarrow Y$ completely contractive maps which 
satisfy conditions 1)-4) for some nets
$\{e_{\alpha}\}\subset X\otimes Y$ and 
$\{f_{\beta}\}\subset Y\otimes X$ with $\|e_{\alpha}\|_{cb} < 1$ and $\|f_{\beta}\|_{cb} < 1$.
 We devide the proof in several steps.

\noindent {\bf Step 1.}
Let $A_1 = X\otimes_h Y$. For each $a\in A_1$ let $L_a : X\rightarrow X$ be the map given by 
$L_a (x) = p(a\otimes x)$. It is clear that $L_a$ is completely contractive. 
Let $A_0=\{a\in A_1 : L_a = 0\}$; it is obvious that $A_0$ is a closed subspace of $A_1$.
We define a multiplication $m$ on $A_1$ as follows. If $a\in A_1$ is an arbitrary element and 
$b = \sum_i x_i\otimes y_i\in A_1$ is an finite sum of elementary tensors, 
then we let $m(a, b) = \sum_i p(a\otimes x_i)\otimes y_i$. We first prove associativity: if 
$a, b=\sum_i x_i\otimes y_i, c=\sum_j u_j\otimes v_j\in A_1$, then 
\begin{eqnarray*}
m(a,m(b,c)) & = & m(a, \sum_{i,j} p(x_i\otimes y_i \otimes u_j)\otimes v_j)\\
                   & = & \sum_{i,j} p(a\otimes p(x_i\otimes y_i \otimes u_j)\otimes v_j\\
                   & = & \sum_{i,j} p(p(a\otimes x_i)\otimes y_i \otimes u_j)\otimes v_j\\
                   & = & \sum_j p(m(a,b)\otimes u_j)\otimes v_j = m(m(a,b),c).
\end{eqnarray*}

By the associativity of the Haagerup tensor
product we have that the map $m$ coincides with the map
$$p\otimes \id : (A_1\otimes_h X) \otimes_h Y \longrightarrow X\otimes_h Y.$$
The last map is completely contractive as a tensor product of completely contractive 
maps, so $m$ is completely contractive. It follows that $m$ can be extended to all of 
$A_1\times A_1$ and the extention, which will be denoted again by $m$, is associative.
Condition 3) implies that $\{e_{\alpha}\}$ is an approximate unit for $A_1$. It follows \cite{brs} 
that $A_1$ is an (abstract) operator algebra. 
Next we point out that $A_0$ is a (two-sided) ideal 
in $A_1$. First observe that for each $a,b\in A_1$ and $z\in X$ we have 
\begin{equation}\label{eq_2}
p(m(a,b)\otimes z) = p(a\otimes p(b\otimes z)).
\end{equation}
Indeed, if $a = x_1\otimes y_1$ and $b=x_2\otimes y_2$ are elementary tensors, then by 1) we have
\begin{eqnarray*}
p(m(a,b)\otimes z) & = & p(p(x_1\otimes y_1\otimes x_2)\otimes y_2\otimes z)\\
                              & = & p(x_1\otimes y_1\otimes p(x_2\otimes y_2\otimes z))=
                                       p(a\otimes p(b\otimes z))
\end{eqnarray*}
and the general case follows by (bi)linearity and approximation.
Now if $a\in A_0$, $x\otimes y\in A_1$ and $z\in X$, then 
$$L_{m(a,x\otimes y)}(z) = L_{p(a\otimes x)\otimes y}(z)=p(p(a\otimes x)\otimes y\otimes z) = 0$$
and, using (\ref{eq_2}),
$$L_{m(x\otimes y, a)}(z) = p(m(x\otimes y, a)\otimes z) = p(x\otimes y\otimes p(a\otimes z)) = 0.$$
This shows that $A_0$ is an ideal in $A_1$. 
Let $A=A_1/A_0$; by \cite{brs} we have that $A$ is an (abstract) operator algebra.
Similarly we define a multiplication on $B_1 = Y\otimes_h X$ which turns $B_1$ into 
an (abstract) operator algebra, 
completely bounded maps $R_b:Y\rightarrow Y$, $b\in B_1$,
an ideal $B_0\subset B_1$ and let $B = B_1/B_0$ be the 
quotient operator algebra.

\noindent {\bf Step 2.}
We equip the spaces $X$ and $Y$ with structures of $A, B$- and $B, A$-bimodules
respectively. If $a\in A_1$ and $x\in X$, let 
$$(a+A_0)\cdot x = p(a\otimes x).$$
If $b\in A_1$ is such that $a-b\in A_0$, then by the definition of $A_0$ we have that 
$p((a-b)\otimes x) = 0$, which implies that the left action of $A$ on $X$ is well defined. From the 
fact that $p$ is completely contractive it follows that this action is completely contractive as well.
The right action of $B$ on $X$ is defined by setting, for each $b\in B$ and $x\in X$,
$$x\cdot (b+B_0) = p(x\otimes b).$$ 
To prove that this is a well defined action, suppose that $b\in B_0$, that is, 
$q(b\otimes y) = 0$ for each $y\in Y$. If $x_1\in X$, $y_1\in Y$, we have that 
$$p(p(x\otimes b)\otimes y_1\otimes x_1) = p(x\otimes q(b\otimes y_1)\otimes x_1) = 0$$
and thus $p(p(x\otimes b)\otimes d)=0$ for each $d\in Y\otimes_h X$ which, according to
4), implies that $p(x\otimes b)=0$. It follows that the right action of $B$ on $X$ is well defined and,
as before, completely contractive. In other words, $X$ is endowed with a structure of an 
$A, B$-operator bimodule. 
The $A, B$-operator bimodule structure of the operator space $Y$ is defined by letting 
$$(b+B_0)\cdot y = q(b\otimes y)$$
and
$$y\cdot (a+A_0) = q(y\otimes a),$$
for each $a\in A_1$, $b\in B_1$ and $y\in Y$. By symmetry, these actions are well-defined
and completely contractive.

\noindent {\bf Step 3.}
We define pairings $(\cdot,\cdot) : X\times Y \rightarrow A$ and 
$[\cdot,\cdot] : Y\times X \rightarrow B$ by letting $(x,y) = x\otimes y +A_0$ and 
$[y,x] = y\otimes x +B_0$. It is obvious that the pairings $(\cdot,\cdot)$ and $[\cdot,\cdot]$
are completely contractive.

\noindent {\bf Step 4. }
We show that the pairings defined in Step 3 are balanced bimodule maps.
Indeed, if $b\in B_1$, $x\in X$ and $y\in Y$, then by 1) we have that 
$p(p(x\otimes b)\otimes y \otimes z) = p(x\otimes q(b\otimes y)\otimes z)$ for each 
$z\in X$, which means that 
$p(x\otimes b)\otimes y - x\otimes q(b\otimes y) \in A_0$ and thus $(x\cdot (b+B_0), y) = 
(x, (b+B_0)\cdot y)$, for each $x\in X$ and $y\in Y$.
Similarly one checks that the pairing $[\cdot,\cdot]$ is $A$-balanced.

Let $a_1= x_1\otimes y_1, a_2 = x_2\otimes y_2\in A_1$ and $z\in X$. Then by 1) we have 
$$p(p(a_1\otimes x)\otimes q(y,a_2)\otimes z) = 
p(p(p(a_1\otimes x)\otimes y\otimes x_2) \otimes y_2 \otimes z)$$
which means that 
$$p(a_1\otimes x)\otimes q(y\otimes a_2) - p(p(a_1\otimes x)\otimes y\otimes x_2) \otimes y_2\in A_0$$
for each $x\in X$ and $y\in Y$ and so $(a_1\cdot x, y\cdot a_2) = a_1 (x,y) a_2$. Similarly, 
the pairing $[\cdot,\cdot]$ is a $B$-bimodule map.

\noindent {\bf  Step 5. }
We have that $(x_1,y)\cdot x_2 = x_1\cdot [y,x_2]$ and 
$[y_1,x]\cdot y_2 = y_1\cdot (x,y_2)$ for each $x_1,x_2\in X$ and $y_1,y_2\in Y$.
These identities are direct consequences of the definitions of the pairings and the module 
actions.

\noindent {\bf Step 6.}
The projection maps $\pi : A_1\rightarrow A$ and $\rho: B_1\rightarrow B$ are complete quotients.
This is an immediate consequence of conditions 3) and 4) and Lemma 2.9. of \cite{bmp}.

Steps 1-6 above ensure that $(A,B,X,Y,(\cdot,\cdot),[\cdot,\cdot])$ is a Morita context with 
completely contractive pairings. 

For the converse direction, suppose that $(A,B,X,Y,(\cdot,\cdot),[\cdot,\cdot])$ is a Morita context with 
completely contractive pairings. Define 
$\tilde{p}:X\times Y\times X\rightarrow X$ and 
$\tilde{q}: Y\times X\times Y\rightarrow Y$ by letting 
$$\tilde{p}(x_1, y, x_2) = (x_1,y)\cdot x_2$$
and 
$$\tilde{q}(y_1, x, y_2) = (y_1,x)\cdot y_2.$$
It is obvious that $\tilde{p}$ and $\tilde{q}$ are completely contractive trilinear maps. If 
$p:X\otimes_h Y\otimes_h X\rightarrow X$ and 
$q: Y\otimes_h X\otimes_h Y\rightarrow Y$ be their linearizations through the Haagerup tensor 
product, then it is easily seen from the definition of a Morita context that 
Properties 1), 2)  3) and 4) hold.

\bigskip

\sproof
Assume that $X$ and $Y$ are operator spaces, 
$p: X\otimes_h Y\otimes_h X\rightarrow X$ and
$q: Y\otimes_h X\otimes_h Y\rightarrow Y$ completely contractive maps which 
satisfy conditions 1)-4) for some families 
$\{e_{\alpha}\}\subset X\otimes Y$ and 
$\{f_{\beta}\}\subset Y\otimes X$ with $\|e_{\alpha}\|_{cb} < 1$ and $\|f_{\beta}\|_{cb} < 1$.
We recall that 
$L_a : X\rightarrow X$ and $R_b : Y\rightarrow Y$ are the operators given by 
$L_a(x) = p(a\otimes x)$ and $R_b(y) = q(b\otimes y)$, where $a\in X\otimes_h Y$ and 
$b\in Y\otimes_h X$. Recall that $L_a\in B(X)$ and $R_b\in B(Y)$. 
Let $\ccl = \{L_a : a\in X\otimes_h Y\}^{-}$ and $\ccr = \{R_b : b\in Y\otimes_h X\}^{-}$, the 
closures being taken in the norm topology of $B(X)$.
As in the First proof, we will only point out the \lq\lq left'' considerations in constructing 
the Morita  context whenever the \lq\lq right'' ones follow either by symmetry or in a similar way.
First note that $\ccl$ is a subalgebra of $B(X)$. Indeed, we have that 
\begin{eqnarray*}
L_{x_1\otimes y_1}L_{x_1\otimes y_1}(z) & = & 
p(x_1\otimes y_1\otimes p(x_2\otimes y_2\otimes z)\\
& = & p(p(x_1\otimes y_1\otimes x_2)\otimes y_2\otimes z) =
L_{p(x_1\otimes y_1\otimes x_2)\otimes y_2} (z)
\end{eqnarray*}
for each $z\in X$ and it follows that $L_aL_b\in \ccl$ whenever $a,b\in X\otimes_h Y$.
Next observe that $\ccl$ is appriximately unital according to 3). Indeed, if $z\in X$, 
$\|z\|\leq 1$ and $a = x\otimes y\in X\otimes_h Y$, then 
\begin{eqnarray*}
\|L_{e_{\alpha}}L_a(z) - L_a(z)\| & = & \|p(p(e_{\alpha}\otimes x)\otimes y\otimes z) -
p(x\otimes y\otimes z)\|\\             & = & \|p((p(e_{\alpha}\otimes x)-x)\otimes y\otimes z)\|\leq
\|p(e_{\alpha}\otimes x)-x\|\|y\|
\end{eqnarray*}
and clearly $\|p(e_{\alpha}\otimes x)-x\|\|y\|\longrightarrow_{\alpha} 0$.

Endow the operator space $X$ with the structure of an $\ccl, \ccr$-bimodule, by letting 
$L_a\cdot x = p(a\otimes x)$ and $x\cdot R_b = p(x\otimes b)$. 
In a similar way, endow the operator space $Y$ with the structure of an $\ccr, \ccl$-bimodule.

Let $n\in\mathbb{N}$. We define a norm $\|\cdot\|_n$ 
on $M_n(\ccl)$, letting, for an element $L = (L_{a_{i j}})\in M_n(\ccl)$, $\|L\|_n$ 
to be the maximum of the numbers
$$\sup \{\|L\odot x\|_n 
: x = (x_{i j})\in M_n(X), \|x\|_n\leq 1\}$$
and
$$\sup \{\|y\odot L\|_n 
: y = (y_{i j})\in M_n(Y), \|y\|_n\leq 1\};$$
it is clear that $\|L\|_n\leq \|p\|_n$.
It is easy to observe that the sequence $\{\|\cdot\|\}_{n\in\mathbb{N}}$ defines an operator 
space structure on $\ccl$. We show that $\ccl$ is actually an operator algebra with respect to 
this operator space structure by 
checking that the multiplication on $\ccl$ is completely contractive. Indeed, if 
$(T_{i j}), (S_{i j})\in M_n(\ccl)$ then 
$$\| (T_{i j})\odot  (S_{i j})\odot (x_{i j})\|_n = \|(\sum_{k,l} T_{i k} (S_{k l} (x_{l j})))_{i, j}\|_n
\leq \| (T_{i j})\|_n\| (S_{i j})\|_n\|(x_{i j})\|_n$$
and, similarly,
$\|((y_{i j})\odot (T_{i j})\odot  (S_{i j})\|\leq 
\| (T_{i j})\|_n\| (S_{i j})\|_n\|(y_{i j})\|_n$,
which proves the assertion.
It follows that $\ccl$ is an operator algebra with respect to the defined sequence of matrix norms.
Similarly, $\ccr$ is an operator algebra with respect to an analogous sequence of matrix norms.
From the definition of the matrix norms 
of $\ccl$ it follows that the module actions of $\ccr$ and $\ccl$ on $X$ and $Y$ 
are completely contractive. 

Define pairings $(\cdot,\cdot) : X\times Y \rightarrow \ccl$ and 
$[\cdot,\cdot] : Y\times X \rightarrow \ccr$ by setting 
$(x,y) = L_{x\otimes y}$ and $[y,x] = R_{y\otimes x}$, where $x\in X$ and $y\in Y$.
The fact that the pairings are ballanced bimodule maps and the \lq\lq linking'' properties 
between the two pairings are verified readily as in the First proof.

The converse direction is established as in the First proof.
\prend

In the case one (both) of the algebras in the Morita context is (are) unital, 
some changes in conditions 
3) or (and) 4) in Theorem \ref{th_mor_pair} must be made. 
For example, suppose that we want to characterize those pairs $(X,Y)$ of operator spaces 
for which there exist operator algebras $A$ and $B$ such that $A$ is unital, $B$ is 
approximately unital and $A$ and $B$ are Morita equivalent with an equivalence pair
$(X,Y)$. Then, from Lemma 2.9 in \cite{bmp} it follows that 
the existence of the family $\{e_{\alpha}\}$ and
condition 3) must be replaced by the following condition:

3') For each $\epsilon > 0$ there exists an element $e_{\epsilon}\in X\otimes_h Y$ such that 
$\|e_{\epsilon}\|<1+\epsilon$ and $p(e_{\epsilon}\otimes x) = x$ for each $x\in X$.

Similar changes must be made in the characterizations of Morita pairs $(X, Y)$ 
between some algebras $A$ and $B$ when both $A$ and $B$ are unital or 
$A$ is approximately unital and $B$ is unital.

Theorem \ref{th_mor_pair} has the following immediate consequence.

\begin{corollary}\label{c_mult}
Let $X$ and $Y$ be operator spaces with the properties 1)-4) of Theorem \ref{th_mor_pair}.
Then there exist completely isometric representations $\nph:X\rightarrow \ccb(\cch,\cck)$ and 
$\psi: Y\rightarrow \ccb(\cck,\cch)$ (where $\cch$ and $\cck$ are Hilbert spaces) such that
$\nph(X)\psi(Y)\nph(X)\subseteq\nph(X)$ and $\psi(Y)\nph(X)\psi(Y)\subseteq\psi(Y)$.
\end{corollary}

\section{The selfadjoint case}

In this section we consider the consequences of Theorem \ref{th_mor_pair} when the 
operator algebras in the Morita context are C*-algebras. If $(X, \{\|\cdot\|_n\}$ 
is an operator space, then
we are able to form its conjugate space $X^*$. As a group with respect to the addition, $X^*$ 
coincides with $X$, the scalar multiplication is given by $\lambda x^* = (\overline{\lambda}x)^*$ 
(where $x^*$ is $x$ viewed as an element of $X^*$). The operator space structure on $X^*$ 
is given by assigning to an $n$ by $n$ matrix $(x^*_{i j})_{i,j}$ with entries in $X^*$ the norm 
$\|(x_{j i})_{i,j}\|_n$. Note that, with respect to these norms, there is a natural identification 
between $M_n(X^*)$ and $M_n(X)^*$.

Recall that an imprimitivity bimodule between two C*-algebras $A$ and $B$ 
is a Banach space $X$ which is an $A,B$-bimodule and
which is equipped with sesquilinear forms
$_A\langle\cdot,\cdot\rangle : X\times X \rightarrow A$ and 
$\langle\cdot,\cdot\rangle_B: X\times X\rightarrow B$, such that $_A\langle\cdot,\cdot\rangle$
is conjugate linear on the second virable, $\langle\cdot,\cdot\rangle_B$ is conjugate linear on the 
first, $(A,X,_A\langle\cdot,\cdot\rangle)$ and $(B,X,\langle\cdot,\cdot\rangle_B)$ are 
Hilbert $A$-modules, $_A\langle x,y\rangle z = x\langle y,z\rangle_B$ and the linear spans of 
$\{_A\langle x,y\rangle : x,y\in X\}$ and $\{\langle x,y\rangle_B : x,y\in X\}$ are dense in $A$ and 
$B$ respectively. If $A$ and $B$ are C*-algebras and $X$ is an imprimitivity bimodule bewteen 
$A$ and $B$, then we are able to form an imprimitivity $B,A$-bimodule $X^*$ in the natural way,
letting $_B\langle x^*,y^*\rangle = \langle y,x\rangle_B^*$ and 
$\langle x^*,y^*\rangle_A = _A\langle y,x\rangle^*$. The set
$$\ccc = 
\left(
\matrix{%
A & X\cr
X^* & B }
\right) $$
can be endowed with a norm turning it into
a C*-algebra under the multiplication performed using the respective module actions and pairings.
Since a C*-algbera has a canonical operator space structure, the inclusion 
$$x\hookrightarrow 
 \left(
\matrix{%
0 & x\cr
0 & 0 }
\right) $$
endows $X$ with a canonical operator space structure. We consider an imprimitivity 
bimodule as an operator space with respect to this canonical operator space structure.

Now we are ready to establish the following theorem.

\begin{theorem}\label{th_C*}
Let $X$ be an operator space. Then $X$ is (completely isometrically isomorphic to) 
an imprimitivity bimodule between certain 
strongly Morita equivalent (in the sense of Rieffel) C*-algebras if and only if there exists a 
completely contractive map
$p : X\otimes_h X^*\otimes_h X\rightarrow X$ such that 

a) $p(p(x_1\otimes x_2^*\otimes x_3)\otimes x_4^*\otimes x_5) = 
p(x_1\otimes p(x_4\otimes x_3^*\otimes x_2)^*\otimes x_5) = 
p(x_1\otimes x_2^*\otimes p(x_3\otimes x_4^*\otimes x_5))$ and

b) $\|p_n ((x_{i j})\odot (x_{i j})^*\odot (x_{i j}))\|_n = \|(x_{i j})\|_n^3$, for each $n\in\mathbb{N}$.
\end{theorem}
\proof
Suppose that an operator space $X$ and a map $p$ are given such that a) and b) are 
satisfied. Let $q:X^*\otimes_h X\otimes_h X^*\rightarrow X^*$ be the map given by 
$$q(x^*\otimes y\otimes z^*) = p(z\otimes y^*\otimes x)^*, \ \ x,y,z\in X.$$
An immediate verification shows that $p$ and $q$ satisfy conditions 1) and 2) of Theorem 
\ref{th_mor_pair}.
As in the Second proof of the proof of Theorem \ref{th_mor_pair} we construct 
the operator spaces $\ccl$ and $\ccr$, which are also algebras with a completely 
contractive multiplication. Note that we are not able to conclude directly that $\ccl$ and $\ccr$
are operator algebras since they are not a priori approximately unital or unital. 
We will show, however, that $\ccl$ and $\ccr$ are C*-algebras and that
their operator space structure, whose sequence of matrix norms will be denoted by 
$\{\|\cdot\|_n\}$, coincides with their canonical C*-algebra operator space structure.
Actually, we will directly show that $M_n(\ccl)$ (and so $M_n(\ccr)$ as well) is a C*-algebra
with respect to $\|\cdot\|_n$. 
First we define involutions on $\ccl$ and $\ccr$:
if $x,y\in X$, we let $L_{x\otimes y^*}^* = L_{y\otimes x^*}$ and 
$R_{y^*\otimes x}^* = R_{x^*\otimes y}$. 
The fact that the mappings $L\rightarrow L^*$ and $R\rightarrow R^*$ are 
indeed involutions on $\ccl$ and $\ccr$ follows easily from condition a).

Let $a_{i j} = \sum_{u=1}^{v_{i j}} x_{i j}^u\otimes y_{i j}^{u *}$,
$L = (L_{a_{i j}})_{i,j}\in M_n(\ccl)$ and $z = (z_{i j})_{i,j}\in M_n(X)$.
Then by b) we have 
\begin{eqnarray*}
& \| & L \ \odot \  z\ \ \ \|_n^3\\ 
& = &
\|p_n \left( ((L_{a_{i j}})\odot (z_{i j}))\odot ((L_{a_{i j}})\odot (z_{i j}))^* \odot 
((L_{a_{i j}})\odot (z_{i j}))\right)\|\\
& = & 
\|p_n ( (\sum_k p(a_{i k}\otimes z_{k j}))_{i,j}\odot (\sum_l p(a_{t l}\otimes z_{l s})^*)_{s, t}
\odot (\sum_m p(a_{r m}\otimes z_{m v}))_{r, v} )\|\\
& = & 
\|\left( (\sum_{j,l,k,r,s,u} p( p(a_{i k}\otimes z_{k j})\otimes  p(a_{s l}\otimes z_{l j})^*
\otimes p(a_{s r}\otimes z_{r v}))\right)_{i, v}\|_n\\
& = & 
\|\left( \sum_{j,l,k,r,s,u} p(p(a_{i k}\otimes z_{k j})\otimes 
p(y_{s r}^u\otimes x_{s r}^{u *}\otimes p(a_{s l}\otimes z_{l j}))^*\otimes z_{r v}\right)_{i, v}\|_n.
\end{eqnarray*}
On the other hand, 
\begin{eqnarray*}
L^*\odot (L\odot z)
& = & 
(L_{a_{j i}}^*)_{i,j}\odot (\sum_l p(a_{s l}\otimes z_{l t}))_{s, t}\\
& = & 
\left( \sum_s L_{a_{s r}}^*(\sum_l p(a_{s l}\otimes z_{l t}))\right)_{r, t}\\
& = & 
\left( \sum_{s,l,u} p(y_{s r}^u\otimes x_{s r}^{u *}\otimes p(a_{s l}\otimes z_{l t}))\right)_{r, t}
\end{eqnarray*}
and thus
\begin{eqnarray*}
p_n ( ( L \odot z ) & \odot & (L^*\odot (L\odot (z)))^{*} \odot z )\\
& = &
p_n\left( \sum_{j,l,k,r,s,u} p(a_{i k}\otimes z_{k j})\otimes p( y_{s r}^u\otimes x_{s r}^{u *}\otimes
p(a_{s l}\otimes z_{l j}))^*\otimes z_{r v}\right)_{i, v}.
\end{eqnarray*}
We conclude that 
\begin{eqnarray*}
\|L\odot z\|_n^3 
& = &
\|p_n ( (L\odot z)  \odot  (L^*\odot (L\odot (z)))^*\odot z )\|_n\\
& \leq & 
\|L\odot z\|_n \|(L^*\odot (L\odot (z))^*\|_n\|z\|_n.
\end{eqnarray*}
Since the involution is isometric, it follows that 
$$\|L\odot z\|^2\leq \|(L^*\odot L)\odot z\|.$$
Similarly, we obtain that 
$$\|z'\odot L\|^2\leq \|z'\odot (L^*\odot L)\|$$
for each $z'\in M_n(Y)$. It follows that 
$$\|L\|_n^2\leq \|L^*\odot L\|_n$$
which implies that $\|\cdot\|_n$ is a C*-norm on $M_n(\ccl)$.
For $n=1$, we obtain that $\ccl$ is a C*-algebra.
Since $M_n(\ccl)$ has a unique C*-norm, it follows that $\|\cdot\|_n$
coincides with the norm on $M_n(\ccl)$ \lq\lq inherited'' from the C*-structure of $\ccl$.

Thus we have shown that $(\ccl,\ccr,X,Y,(\cdot,\cdot),[\cdot,\cdot])$ is a Morita context. 
From Theorem 6.2 of \cite{bmp} it follows that there exists a complete isometry 
$i: X^*\longrightarrow X^*$ such that $X$ becomes an imprimitivity bimodule between $\ccl$
and $\ccr$ under the pairings $_{\ccl}\langle x,y\rangle = (x,i(y^*))$ and 
$\langle x,y\rangle_{\ccr} = [i(x^*),y]$, $x,y\in X$. It follows that, as a Banach space, 
$X$ is an imprimitivity bimodule in the sense of Rieffel. 
It only remains to be checked that the operator space structure of $X$ coincides with the 
canonical imprimitivity bimodule operator space structure $\{\|\cdot\|'_n\}$ on $X$. 
By Lemma 5.6 of \cite{bmp} we have that there are completely isometric representations
$\sigma,\pi,\nph,\psi$ of $\ccl,\ccr, X, X^*$ on appropriate Hilbert spaces 
such that if $\rho$ is the respective two by two matrix 
representation formed by $\sigma,\pi,\nph$ and $\psi$, then $\rho$ is completely isometric
and $\sigma ((x,y^*)) = \nph(x)\psi(y^*)$ for each $x,y\in X$.
Let $x\in X$. Then
\begin{eqnarray*}
\|x\|^{' 2}_1 & = & \|_{\ccl}\langle x,i(x^*)\rangle\| = \|\sigma ((x,i(x^*)))\|\\
                   & = & \|\nph(x)\psi(i(x^*))\| = \|\nph(x)\nph(x)^*\| = \|\nph(x)\|^2 = \|x\|^2.
\end{eqnarray*}
Thus, $\|\cdot\|'=\|\cdot\|$. Note that for each $n\in\mathbb{N}$ the pair
$(M_n(X),M_n(X^*))$ is a Morita equivalence pair between the C*-algebras $M_n(\ccl)$ and 
$M_n(\ccr)$. Just as above we conclude that $\|\cdot\|^{'}_n=\|\cdot\|_n$.

The converse direction is obtained as the converse direction in Theorem 
\ref{th_mor_pair}.
\prend

Note that there was no need to impose a condition about existence of approximate 
units in Theorem \ref{th_C*}: as we saw, the non-degeneracy condition a) assures that 
the algebras we construct are C*-algebras; a posteriori they are approximately unital or unital.

\smallskip

Theorem \ref{th_C*} leads us to the following definition.

\begin{definition}\label{d_tro}
A ternary operator system is 
an operator space $X$ equipped with a triple product $[\cdot,\cdot,\cdot] : X\times X\times X
\longrightarrow X$, linear on the first and the third variable and conjugate linear on the second
such that 

a) $[[x_1,x_2,x_3],x_4,x_5] = [x_1,[x_4,x_3,x_2],x_5] = [x_1,x_2,[x_3,x_4,x_5]]$ and

b) $\|[(x_{i j})\odot (x_{ j i}) \odot (x_{i j})]_n\| = \|(x_{i j})\|_n^3$, for each  $(x_{i j})\in M_n(X)$, 
$n\in\mathbb{N}$.
\end{definition}

Recall that a ternary C*-ring is a Banach space endowed with a contractive triple product
$[\cdot,\cdot,\cdot] : X\times X\times X\longrightarrow X$,
linear on the first and the third variable and conjugate linear on the second, satisfying 
conditions a) and  b) in Definition \ref{d_tro} for $n=1$. 
A consequence of Zettl's results is that the Banach spaces which are imprimitivity bimodules
between certain C*-algebras are, up to isometry, 
precisely the ternary C*-rings. But, as we pointed out, an imprimitivity bimodule has a 
canonocal structure of an operator space.
Theorem \ref{th_C*} gives an operator space version of the above result: it states 
that the operator spaces which are imprimitivity bimodules 
bewteen certain C*-algebras 
are, up to a complete isometry, precisely the ternary operator systems. The conditions 
a) and b) for n=1 are not sufficient to ensure that the operator space is {\it completely}
isometric to an imprimitivity bimodule - we must require moreover that condition b) is
fulfilled for every $n\in\mathbb{N}$. An example of an
operator space $X$ which is a ternary C*-ring but not a ternary operator system
is a C*-algebra equipped 
with some non-canonical operator space structure.

Ternary C*-rings and ternary operator systems are the linear space analogue of the 
C*-algebras. A concrete version of these objects are the so called {\it ternary rings of operators
(TRO)}. Recall \cite{z} that a TRO is a 
subspace $V\subseteq \ccb(H, K)$ of operators between 
some Hilbert spaces $H$ and $K$, closed under the (canonical) triple product 
$(a,b,c)\longrightarrow ab^*c$. 
Note that, as a concrete operator space, a TRO has a natural operator space structure. 
It is clear that a TRO is a ternary C*-ring and a ternary operator 
system with respect to the canonical triple product.
A result analogous to the Gelfand-Naimark Theorem is proved in \cite{z}: Each 
ternary C*-ring $X$ has a decomposition $X=X_+\oplus X_-$ such that $X_+$ is 
(isometrically) isomorphic to 
a TRO while $X_-$ is (isometrically) anti-isomorphic to a TRO. 
Of course, a morphism between ternary C*-rings is defined as a bounded linear map which preserves
the triple product. Similarly, a morphism between ternary operator systems is a 
completely bounded linear map which preserves the triple product.
Theorem \ref{th_C*} has the following consequences.


\begin{corollary}\label{c_tro}
Let $X$ be a ternary operator system. Then there exist ternary operator sub-systems 
$X_+$ and $X_-$ such that $X=X_+\oplus X_-$ and $X_+$ is completely isometrically 
isomorphic to a TRO while $X_-$ is completely isometrically anti-isomorphic to a TRO.
\end{corollary}

Thus one is able to decompose a ternary operator system into 
parts on which the triple product can be given a concrete expression through Hilbert space 
operators operations. 
The triple product on the space $X_+$ above is given, up to complete isometry, by 
$(a,b,c)\longrightarrow ab^*c$ while the triple product on $X_-$ is given, again up to 
complete isometry, by $(a,b,c)\longrightarrow -ab^*c$. Theorem \ref{th_C*} shows that, 
up to complete isometry, imprimitivity bimodules do not recognize the difference between 
these two basic triple products.

If we are not interested in preserving the triple products, we can achieve the following
representation.

\begin{corollary}\label{c_tro1}
Let $X$ be a ternary operator system. Then $X$ is completely isometric to a TRO.
\end{corollary}

\bigskip

\noindent{\bf Acknowledgement }
The author wishes to thank Aristides Katavolos for the useful discussions on the topic of this paper.

\bigskip



\begin{thebibliography}{99}

\bibitem{brs} D. P. Blecher, Z-j Ruan and A.M. Sinclair, A Characterization of Operator Algebras,
{\it J. Funct. Anal.}, 89 (1990), 188-201

\bibitem{bmp} D. P. Blecher, P. S. Muhly and V. I. Paulsen, Categories of Operator Modules
(Morita Equivalence and Projective Modules), {\it Mem. Amer. Math. Soc.}, 143 (2000), no.681

\bibitem{cs} E. Christensen and A. M. Sinclair, Representations of completely bounded 
multilinear operators, {\it J. Functional Anal.} 72 (1987), 151-181

\bibitem{er} E. Effros, Z-j. Ruan, {\it Operator Spaces}, Clarendon Press, Oxford, 2000


\bibitem{rief} M. Rieffel, Morita equivalence for C*-algebras and W*-algebras, 
{\it J. Pure Appl. Algebra} 5 (1974), 51-96

\bibitem{z} H. Zettl, A Characterization of Ternary Rings of Operators, 
{\it Adv. Math}, 48 (1983), 117-143

\end{thebibliography}
\end{document}